# Important Notes on Lyapunov Exponents

Keying Guan

Sciense College, Beijing Jiaotong University,

Beijing, China, 100044

Email: keying.guan@gmail.com

**Abstract:** It is shown that the famous Lyapunov exponents cannot be used as the numerical characteristic for distinguishing different kinds of attractors, such as the equilibrium point、the limit closed curve、the stable torus and the strange attractor.

## 1. Introduction

The conception "Lyapunov Exponent" has been used widely in the study of dynamical system.

Usually, the **Lyapunov exponent** or **Lyapunov characteristic exponent** of a dynamical system is a quantity that characterizes the rate of separation of infinitesimally close trajectories $\mathbf{Z}(t)$ and $\mathbf{Z_0}(t)$ in phase space.

Let $\delta\mathbf{Z}(t) = \mathbf{Z}(t) - \mathbf{Z_0}(t)$, $\delta\mathbf{Z_0} = \mathbf{Z}(0) - \mathbf{Z_0}(0)$, if

$$|\delta\mathbf{Z}(t)| \approx e^{\lambda t}|\delta\mathbf{Z_0}| \qquad (1)$$

then $\lambda$ is treated as the Lyapunov exponent.

If the trajectory $\mathbf{Z}(t)$ is given by a n-dimensional linear ordinary differential equation system with constant coefficients

$$\dot{\mathbf{Z}} = \mathbf{A}\,\mathbf{Z} + \mathbf{f}(t) \qquad (2)$$

and if the constant coefficient matrix $\mathbf{A}$ has n different eigenvalues, $\lambda_1, \lambda_2, \cdots, \lambda_n$, then the real parts of the n different eigenvalues are naturally the Lyapunov exponents. However, if the dynamical system is not given by (2), for instance, if the dynamical system is a nonlinear polynomial autonomous system, the conception on Lyapunov exponents becomes complicated.

Consider the commonly accepted definition of Lyapunov Exponents, which are quoted from the reference [1], and from the Wikipedia on "Lyapunov exponent" (http://en.wikipedia.org/wiki/Lyapunov_exponent).

A. **The maximal Lyapunov exponent** can be defined as follows:

$$\lambda_{max} = \lim_{t \to \infty} \lim_{\delta \mathbf{Z}_0 \to 0} \frac{1}{t} \ln \frac{|\delta \mathbf{Z}(t)|}{|\delta \mathbf{Z}_0|} \qquad (3)$$

The limit $\delta \mathbf{Z}_0 \to 0$ ensures the validity of the linear approximation at any time. It is required that the two limits cannot be exchanged, otherwise, in bounded attractors, the result would be trivially 0.

B. For a dynamical system with evolution equation $f^t$ in an $n$-dimensional phase space, **the spectrum of Lyapunov exponents**

$$\{\lambda_1, \lambda_2, \cdots, \lambda_n\}$$

in general, depends on the starting point $x_0$.

The Lyapunov exponents describe the behavior of vectors in the tangent space of the phase space and are defined from the Jacobian matrix

$$J^t(x_0) = \left. \frac{df^t(x)}{dx} \right|_{x_0} \qquad (4)$$

The $J^t$ matrix describes how a small change at the point $x_0$ propagates to the final point $f^t(x_0)$. The limit

$$\lim_{t \to \infty} (J^t \cdot \text{Transpose}(J^t))^{1/2t} \qquad (5)$$

defines a matrix $\mathbf{V}(x_0)$ (the conditions for the existence of the limit are given by the Oseledec theorem). If $\Lambda_i(\mathbf{x}_0)$ are the eigenvalues of $\mathbf{V}(x_0)$, then the Lyapunov exponents $\lambda_i$ are defined by

$$\lambda_i(x_0) = \operatorname{Ln} \Lambda_i(x_0) \tag{6}$$

Based on the experience of the linear system (2) and some plausible thinking, for a dissipative system, as criterions, it is proposed in the reference [1] that, if the attractor reduces to

(a) stable fixed point, all the exponents are negative;
(b) limit cycle, an exponent is zero and the remaining ones are all negative;
(c) $k$-dimensional stable torus, the first $k$ LEs vanish and the remaining ones are negative;
(d) for strange attractor generated by a chaotic dynamics at least one exponent is positive.

The above-mentioned definition on Lyapunov exponents and proposed criterions about the relations between the characteristic of LE and the properties of the attractors are widely used.

## 2. Some notes on the definition of Lyapunov Exponents

It should be point out first that the expression (3) is not strictly, since the value of right hand side of it is usually not uniquely determined for a given trajectory. In fact the right hand side of (3) depends on how $\delta \mathbf{Z}_0 \to 0$. To see this fact, without lost of generality, assume that the given dynamical system is a three-dimensional autonomous system

$$\begin{cases} \dfrac{dx}{dt} = P(x, y, z) \\ \dfrac{dy}{dt} = Q(x, y, z) \\ \dfrac{dz}{dt} = R(x, y, z) \end{cases} \tag{7}$$

Or for simple, (7) can be written as

$$\frac{d\mathbf{r}}{dt} = \mathbf{F}(\mathbf{r}) \tag{7'}$$

where $\mathbf{r} = (x, y, z)^T$, $\mathbf{F} = (P, Q, R)^T$

Let

$\mathbf{r}_0(t) = (x_0(t), y_0(t), z_0(t))^T$ and $\mathbf{r}(t) = (x(t), y(t), z(t))^T$

are two trajectories of (7). Formally, they satisfy respectively

$$\mathbf{r}_0(t) = \mathbf{r}_0(0) + \int_0^t \mathbf{F}(\mathbf{r}_0(s))ds \tag{8}$$

$$\mathbf{r}(t) = \mathbf{r}(0) + \int_0^t \mathbf{F}(\mathbf{r}(s))ds \tag{9}$$

Let $\delta\mathbf{r}(t) = \mathbf{r}(t) - \mathbf{r}_0(t) = (\delta x(t), \delta y(t), \delta z(t))$, and

$$\delta\mathbf{r}_0 = \mathbf{r}(0) - \mathbf{r}_0(0) = (\delta x(0), \delta y(0), \delta z(0))^T = \delta l (\cos\alpha, \cos\beta, \cos\gamma)^T$$

(10)

where $\delta l = \sqrt{(\delta x(0))^2 + (\delta y(0))^2 + (\delta z(0))^2}$, and $(\cos\alpha, \cos\beta, \cos\gamma)^T$ is the unit vector of $\delta\mathbf{r}_0$ represented with its azimuth. When $\delta l$ is small enough, that is, these two trajectories are close enough, then, for fixed $t$, $\delta\mathbf{r}(t)$ can be treated a small quantity with the same order as $\delta\mathbf{r}_0$. So

$$\delta\mathbf{r}(t) = \exp(\int_0^t \mathbf{J}(\mathbf{r}_0(s))ds)\, \delta\mathbf{r}_0 + \mathbf{o}(\delta\mathbf{r}_0) \tag{11}$$

where $\mathbf{o}(\delta\mathbf{r}_0)$ is a higher order small quantity of $\delta\mathbf{r}_0$ (ref.[2]).

Therefore, when $|\delta\mathbf{r}_0| \to 0$ along fixed direction $(\cos\alpha, \cos\beta, \cos\gamma)$, then

$$\lim_{|\delta\mathbf{r}_0| \to 0} \frac{|\delta\mathbf{r}(t)|}{|\delta\mathbf{r}_0|} = \left| \exp(\int_0^t \mathbf{J}(\mathbf{r}_0(s))ds)\, (\cos\alpha, \cos\beta, \cos\gamma)^T \right| \tag{12}$$

This means the limit (12) depends obviously on the direction of $|\delta\mathbf{r}_0| \to 0$.

Let $J_1(t), J_2(t), \cdots, J_n(t)$ be the n eigenvalues of the matrix $\int_0^t \mathbf{J}(\mathbf{r}_0(s))ds$, and assume $J^*(t)$ is the max real part of these eigenvalus. Then

$$\max[\lim_{|\delta\mathbf{r}_0| \to 0} \frac{|\delta\mathbf{r}(t)|}{|\delta\mathbf{r}_0|}] = \exp J^*(t) \tag{13}$$

From (13), if the following limit exists

$$\lim_{t \to \infty} \frac{1}{t} \ln \max[\lim_{|\delta\mathbf{r}_0| \to 0} \frac{|\delta\mathbf{r}(t)|}{|\delta\mathbf{r}_0|}] = \lim_{t \to \infty} \frac{J^*(t)}{t} \tag{14}$$

then this limit can be reasonably treated as the maximal Lyapunov exponent $\lambda_{max}$. And the $n$ limits

$$\lambda_i = \lim_{t \to \infty} \frac{\operatorname{Re} J_i(t)}{t} \qquad (15)$$

may called the Lyapunov exponents of the trajectory $\mathbf{r}_0(t)$. This definition has obvious meaning and has a close relation with the expected relation (1).

However, there are some issue should be discussed on the commonly used Lyapunov exponents defined by (5) and (6).

Clearly, $J^t$ is equivalent to $\exp(\int_0^t \mathbf{J}(\mathbf{r}_0(s))ds)$, and

$$(J^t \cdot \operatorname{Transpose}(J^t))^{1/2t} = \exp(\frac{\int_0^t \mathbf{J}(\mathbf{r}_0(s))ds + \operatorname{Transpose}(\int_0^t \mathbf{J}(\mathbf{r}_0(s))ds)}{2t})$$

So by (5) and (6), the Lyapunov exponents are just the eigenvalues of the following limit matrix

$$\lim_{t \to \infty} \frac{\int_0^t \mathbf{J}(\mathbf{r}_0(s))ds + \operatorname{Transpose}(\int_0^t \mathbf{J}(\mathbf{r}_0(s))ds)}{2t} \qquad (16)$$

These exponents may be just the same as those given by (15), if $\int_0^t \mathbf{J}(\mathbf{r}_0(s))ds$ is with some symmetries. For instance, if $\int_0^t \mathbf{J}(\mathbf{r}_0(s))ds = (\int_0^t \mathbf{J}(\mathbf{r}_0(s))ds)^T$, or if

$$\int_0^t \mathbf{J}(\mathbf{r}_0(s))ds = \begin{pmatrix} at & 0 & 0 \\ 0 & bt & ct \\ 0 & -ct & bt \end{pmatrix}$$

But, in most cases, the two kinds of Lyapunov exponents are different. For example, if

$$\int_0^t \mathbf{J}(\mathbf{r}_0(s))ds = \begin{pmatrix} -t & 0 & 0 \\ 0 & -2t & t \\ 0 & -t & -3t \end{pmatrix}$$

the exponents given by (5) and (6) are $-1$, $-2$ and $-3$, and the exponents given by (15) are $-1, -\frac{5}{2}$ and $-\frac{5}{2}$.

Sometimes, the difference between the two kinds of exponents is substantive. For instance, if

$$\int_0^t \mathbf{J}(\mathbf{r}_0(s))ds = \begin{pmatrix} -t & 10\,t & 0 \\ 0 & -2\,t & 0 \\ 0 & 0 & -3t \end{pmatrix}$$

then the exponents given by (5) and (6) are $\dfrac{\sqrt{101}-3}{2}, -3$ and $\dfrac{-\sqrt{101}-3}{2}$, and the exponents given by (15) are $-1, -2$ and $-3$. In the first group of the exponents, the maximum one is $\dfrac{\sqrt{101}-3}{2}$, which is positive, and in the second group, the maximum one is $-1$. The two maximum exponents indicate two opposite stabilities.

Because of the difference mentioned above, the exponents given by (15) will be denoted as $LE_J$, and the exponents given by (5) and (6) will be denoted as $LE_O$.

Clearly, in the application of the Lyapunov exponents, the existence of these numbers is very important. Oseledec proved that the limit matrix $\mathbf{V}(x_0)$ of (5) exists with the exception of a subset of initial conditions of zero measure (ref. [2]). This fact might be the advantage of the definition of the exponents given by (5) and (6).

However, the author believes that $LE_O$ has lost the basic practical significance when it is different to $LE_J$.

In the case that the trajectory $\mathbf{r}_0(t)$ is a closed orbit, that is, $\mathbf{r}_0(t)$ is periodic, the limit

$$\lim_{t \to \infty} \frac{\int_0^t \mathbf{J}(\mathbf{r}_0(s))ds}{t} = \frac{\int_0^T \mathbf{J}(\mathbf{r}_0(s))ds}{T} \qquad (17)$$

does exist. So, both of $LE_J$ and $LE_O$ exist. In the following section, we will study what these exponents can be for the spatial limit closed orbits (including the spatial limit cycle).

## C. The Lyapunov exponents for some spatial limit close orbit

In order to get some exact results, this paper will study first some limit cycles, which can be represented exactly with simple elementary functions.

A spatial limit cycle is said meta-stable in this paper, if there are some trajectories approaching to the cycle, and in its any small neighborhood, there are also some trajectories going away from the neighborhood as $t \to +\infty$.

Consider three-dimensional autonomous system

$$\begin{cases} \dfrac{dx}{dt} = y + x(1 - x^2 - y^2)(1 + \alpha - x^2 - y^2) \\ \dfrac{dy}{dt} = -x + y(1 - x^2 - y^2)(1 + \alpha - x^2 - y^2) \\ \dfrac{dz}{dt} = \beta^2 z - z^3 \end{cases} \quad (18)$$

where $\alpha$ ( $\alpha^2 < 1$ ) and $\beta$ are real parameters. The system (18) has 6 different limit cycles:

(i) $x = \sin t, \ y = \cos t, \ z = 0$

(ii) $x = \sqrt{1+\alpha} \sin t, \ y = \sqrt{1+\alpha} \cos t, \ z = 0$

(iii) $x = \sin t, \ y = \cos t, z = \beta$

(iv) $x = \sin t, \ y = \cos t, z = -\beta$

(v) $x = \sqrt{1+\alpha} \sin t, \ y = \sqrt{1+\alpha} \cos t, \ z = \beta$

(vi) $x = \sqrt{1+\alpha} \sin t, \ y = \sqrt{1+\alpha} \cos t, \ z = -\beta$

Since the system (18) is strongly symmetrical, $LE_J$ and $LE_O$ are just the same.

For the limit cycle (i), the eigenvalues of (17) are $\beta^2, -\alpha + i$ and $-\alpha - i$, and $LE_J = LE_O$: $\beta^2, -\alpha, -\alpha$.

If $\beta \ne 0$, the limit cycle is meta-stable when $\alpha \ge 0$, and it is unstable when $\alpha < 0$. If $\beta = 0$, the limit cycle is asymptotically stable when $\alpha > 0$, and it is metastable when $\alpha \le 0$.

For the limit cycle (ii), the eigenvalues of (17) are $\beta^2, \alpha^2 + \alpha + i, \alpha^2 + \alpha - i$, and $LE_J = LE_O$: $\beta^2, \alpha^2 + \alpha, \alpha^2 + \alpha$.

If $\beta \neq 0$, the limit cycle is unstable when $\alpha > 0$, and it is meta-stable when $\alpha \leq 0$. If $\beta = 0$, the limit cycle is meta-stable when $\alpha \geq 0$, and it is asymptotically stable when $\alpha < 0$.

For the limit cycle *(iii)* and *(iv)*, the eigenvalues of (17) are $-2\beta^2$, $-\alpha + i$ and $-\alpha - i$, and LE$_J$ = LE$_O$: $-2\beta^2, -\alpha, -\alpha$.

The limit cycle is asymptotically stable when $\alpha > 0$, and it is meta-stable when $\alpha \leq 0$.

For the limit cycle *(v)* and *(vi)*, the eigenvalues of (17) are $-2\beta^2$, $\alpha^2 + \alpha + i$ and $\alpha^2 + \alpha - i$, and LE$_J$ = LE$_O$: $-2\beta^2, \alpha^2 + \alpha, \alpha^2 + \alpha$.

The limit cycle is meta-stable when $\alpha \geq 0$, and it is asymptotically stable when $\alpha < 0$.

From the above results, it is easy to see the signs of Lyapunov exponents for a asymptotically stable limit cycle may have two kinds of distribution.

(a') All of the three Lyapunov exponents are negative, i.e.,

$$(-, -, -).$$

It happens in the case (iii) and (iv) when $\beta \neq 0$ and $\alpha > 0$, and also in the case (v) and (vi) when $\beta \neq 0$ and $\alpha < 0$,

(b') One of the Lyapunov exponents is zero, and the other two are negative, i.e.,

$$(0, -, -).$$

It happens in the case (i) when $\beta = 0$ and $\alpha > 0$, in the case (ii) when $\beta = 0$ and $\alpha < 0$, and in the case (v) and (vi) when $\beta = 0$ and $\alpha < 0$.

Besides, there are still some other possible sign distributions of the Lyapunov exponents, such that the corresponding limit cycle is asymptotically stable.

Consider the following system

$$\begin{cases} \dfrac{dx}{dt} = y + x(1-x^2-y^2)^3 \\ \dfrac{dy}{dt} = -x + y(1-x^2-y^2)^3 \\ \dfrac{dz}{dt} = \beta^2 z - z^3 \end{cases} \quad (19)$$

It is easy to see that this system has three limit cycles if $\beta \neq 0$:

(i') $x = \sin t, \ y = \cos t, \ z = 0$

(ii') $x = \sin t, \ y = \cos t, \ z = b$

(iii') $x = \sin t, \ y = \cos t, \ z = -b$

In the case (ii') and (iii') when $\beta \neq 0$, the limit cycles are both asymptotically stable, and $LE_J = LE_O$: $0, 0, -2\beta^2$.

In the case (I') when $\beta = 0$, the limit cycle is also asymptotically stable, and $LE_J = LE_O$: $\{0, 0, 0\}$.

From this example we see that there are other two possible sign distributions of the Lyapunov exponents, such that the corresponding limit cycle is asymptotically stable:

(c') Two zeros and one negative, i.e.,

(0, 0, −)

(d') Three zeros, i.e.,

(0, 0, 0)

Note: From the above examples, it is easy to see that, some exponents equals zero is usually corresponding to the critical situation that the stability of the limit cycle changes, or that the number of asymptotically limit cycles changes. It seems a very reasonable criterion for the determination of the stability of the limit cycle. But we will seem that this criterion may not be true in general.

In all of the above examples, $LE_J$ and $LE_O$ are the same. In the following example, more unexpected facts on the Lyapunov exponents for limit closed orbits will appear.

Consider the particular Silnikov equation system

$$\begin{cases} \dfrac{dx}{dt} = y \\ \dfrac{dy}{dt} = z \\ \dfrac{dz}{dt} = x^3 - a^2 x - y - b\,z \end{cases} \quad (20)$$

where the parameters $a$ and $b$ are positive.

In [3] and [4], this system is proved to be an ideal system for the study of three-dimensional differential dynamical systems since it has different kinds of attractors, including spatial limit closed orbits with different rotation numbers.

It has been shown that when the parameter $b$ is slightly smaller than $a^2$, the system (2) will have a limit cycle around the origin, i.e., the Hopf bifurcation (Beiye Feng[1] and then Hongwei Liu[2] have given respectively the strict proof for the existence of the Hopf bifurcation).

This system is not integrable with quadrature, just as expected by most researchers. In fact, it has been proven by Yanxia Hu[3] that this system does not admit any global analytical Lie group, except the trivial one: $t \to t + c$.

Therefore, there is no hope for us to represent the limit cycle of (20) with the quadrature.

So, the calculation of the Lyapunov exponents for its limit cycles, or more generally for the limit closed orbits with different rotation numbers can only be realized numerically.

The following is a series of numerical results of different limit closed orbits of (20):

($n_1$) $a = 1$ and $b = 0.8$. The (20) has a limit cycle of period $T \cong 6.2848$ (see Fig. 1). For the limit cycle,

---

1) Beiye Feng, 关于在b=1 时发生Hopf 分支的证明, http://blog.sciencenet.cn/blog-553379-750729.html
2) Hongwei Liu, et al., 一类Silnikov方程的Hopf分岔及其稳定性, to appear
3) Yanxia Hu, The Non-integrability of a Silnikov Equation, to appear

LE$_J$: −0.0701, −0.0701, −0.6600  (−,−,−)

LE$_O$: 0.5347, −0.4003, −0.9345  (+,−,−)

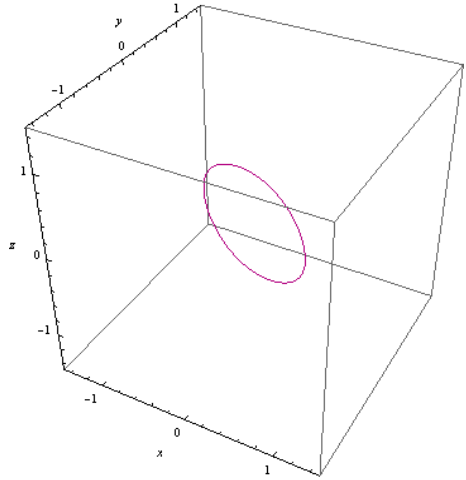

Fig.1   $a = 1$ and $b = 0.8$.

($n_2$) $a = 1$ and $b = 0.6$.  The system (20) has a limit cycle of period $T \cong 6.2899$ (see Fig. 2).  For this limit cycle,

LE$_J$: −0.1929, −0.1929, −0.2142  (−,−,−)

LE$_O$: 0.5044, −0.4654, −0.6390  (+,−,−)

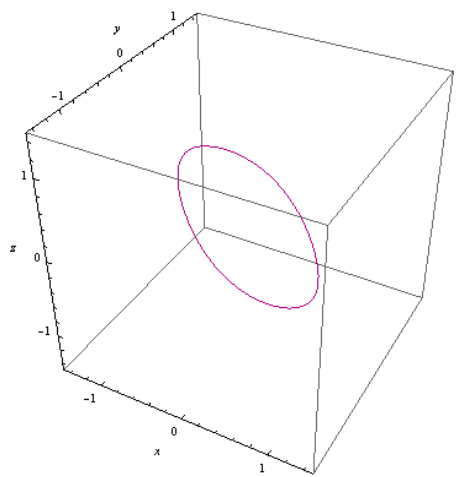

Fig. 2   $a = 1$ and $b = 0.6$

($n_3$) $a = 1$ and $b = 0.5024$.  The system (20) has a limit cycle of period

$T \cong 6.2938$ (see Fig. 3).  For this limit cycle,

    LE$_J$:  $-0.00018,\ -0.2511,\ -0.2511$                                                     $(-,-,-)$

    LE$_O$:  $0.5000,\ -0.5000,\ -0.5024$                                                            $(+,-,-)$

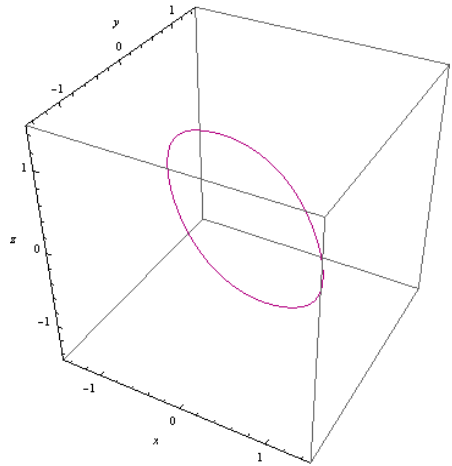

Fig.3    $a = 1$ and $b = 0.6$

($n_4$) $a = 1$ and $b = 0.5023$. The system (20) has still one limit cycle of period $T \cong 6.2938$ (see Fig. 4).  For the limit cycle,

    LE$_J$:  $0.000017,\ -0.2512,\ -0.2512$                                                   $(+,-,-)$

    LE$_O$:  $0.5000,\ -0.5000,\ -0.5023$                                                      $(+,-,-)$

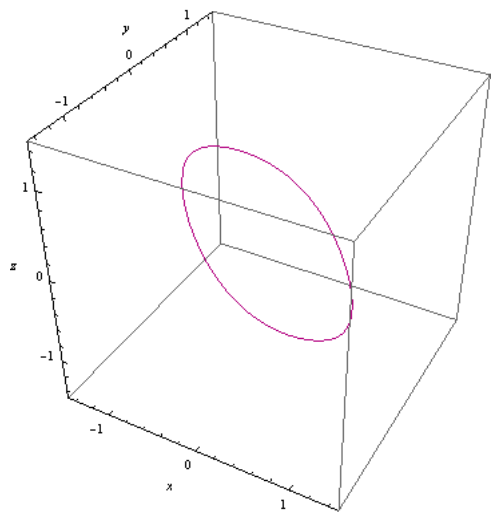

Fig.4    $a = 1$ and $b = 0.5023$.

($n_5$) $a = 1$ and $b = 0.5000$.  The system (20) has still one limit cycle of period $T \cong 6.2939$ (see Fig. 5).  For this limit cycle,

    $LE_J$:   0.0046, $-0.2523$, $-0.2523$                                                    $(+,-,-)$

    $LE_O$:   0.5000, $-0.4984$, $-0.5016$                                                    $(+,-,-)$

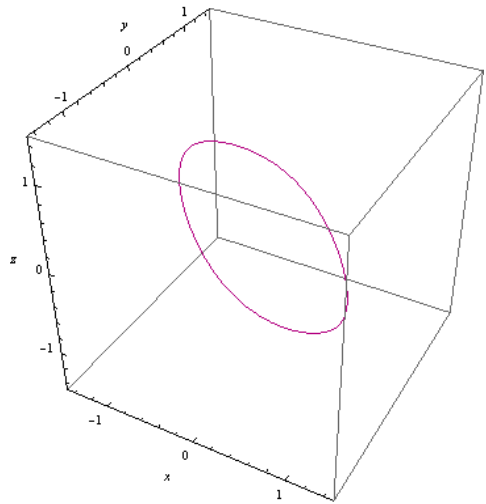

Fig.5    $a = 1$ and $b = 0.5000$

($n_6$) $a = 1$ and $b = 0.4900$.  The system (20) has still one limit cycle of period $T \cong 6.2944$ (see Fig. 6).  For this limit cycle,

    $LE_J$:   0.0244, $-0.2572$, $-0.2572$                                                    $(+,-,-)$

    $LE_O$:   0.5000, $-0.4850$, $-0.5051$                                                    $(+,-,-)$

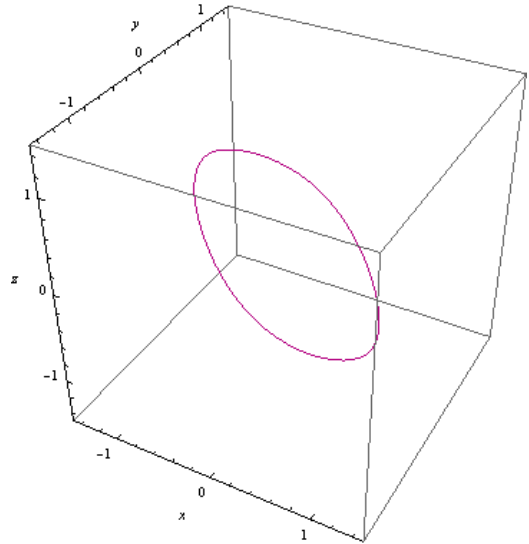

Fig. 6   $a = 1$ and $b = 0.4900$

($n_7$) $a = 1$ and $b = 0.4893$.   Numerically, the system (20) has still a limit cycle of period $T \cong 6.2944$ (see Fig. 7).   For this limit cycle,

$LE_J$:   0.0258, −0.2576, −0.2576         (+,−,−)

$LE_O$:   0.5001, −0.4840, −0.5054         (+,−,−)

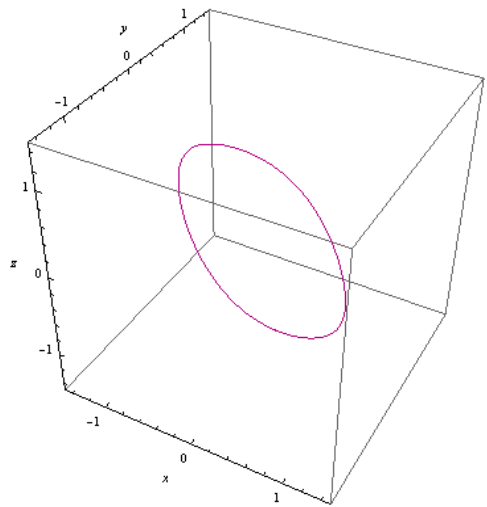

Fig. 7   $a = 1$ and $b = 0.4893$.

($n_8$) $a = 1$ and $b = 0.4892$.   From the numerical result, it can be seen that the number of asymptotically stable limit cycles of system (20) has become two (see

Fig. 8). They are symmetrical about the origin and are very close to each other. The period of each limit cycle is $T \cong 6.2944$. They have the same $LE_J$ and the same $LE_O$,

    $LE_J$:  0.0259, −0.2576, −0.2576                        (+,−,−)

    $LE_O$:  0.5001, −0.4839, −0.5054                        (+,−,−)

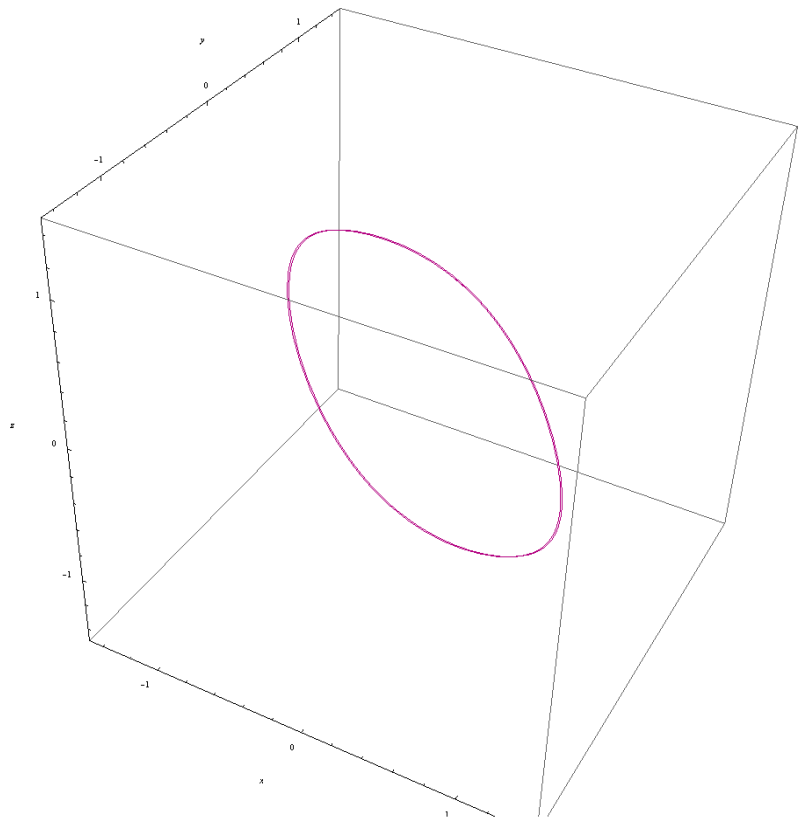

    Fig. 8    $a = 1$ and $b = 0.4892$.   two slightly separated limit cycles

(n$_9$) $a = 1$ and $b = 0.4981$.   The system (20) has two asymptotically stable clearly separated and symmetrical limit cycles (see Fig. 9).   The period of each limit cycle is $T \cong 6.2945$.   For them

    $LE_J$:  0.0260, −0.2576, −0.2576                        (+,−,−)

    $LE_O$:  0.5001, −0.4838, −0.5054                        (+,−,−)

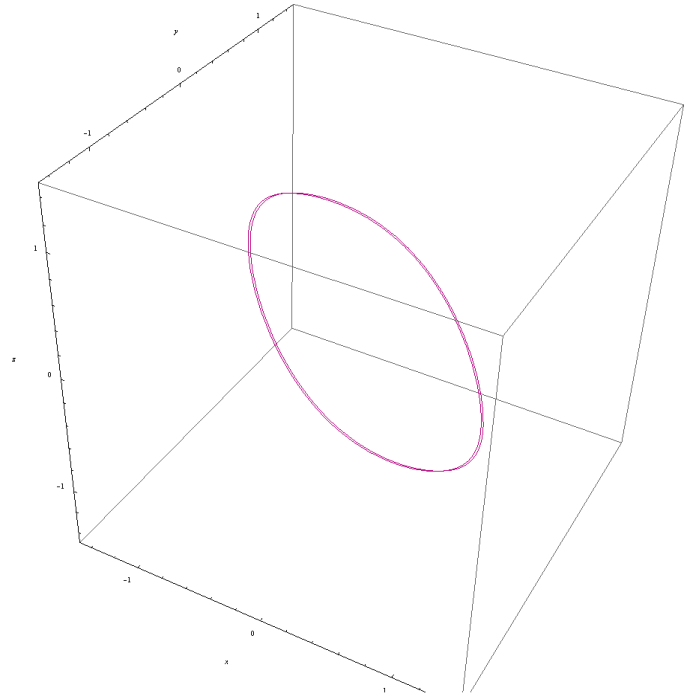

Fig. 9  $a = 1$ and $b = 0.4981$, two clearly separate limit ccycles

($n_{10}$) $a = 1$ and $b = 0.3920$.  The system (20) has two asymptotically stable clearly separated and symmetrical limit closed orbits. Their rotation number are both two (see Fig. 10). The period of each limit cycle is $T \cong 12.7176$.  They have the same LE$_J$ and the same LE$_O$,

LE$_J$:  0.1086, −0.2503, −0.2503 $\qquad\qquad\qquad$ (+,−,−)

LE$_O$:  0.5018, −0.3802, −0.5137 $\qquad\qquad\qquad$ (+,−,−)

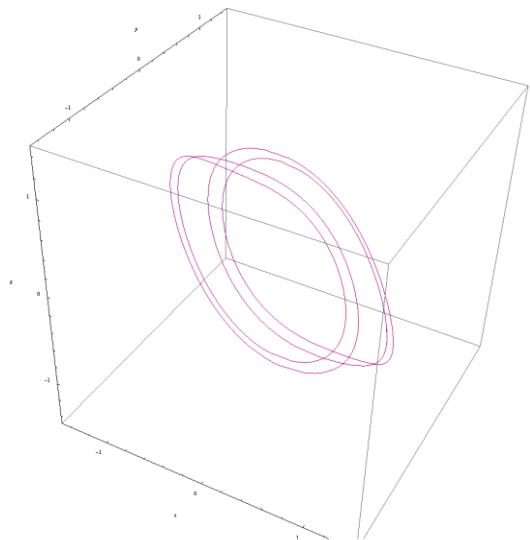

Fig. 10  $a = 1$ and $b = 0.3920$ , two limit closed orbits with rotation number 2

($n_{11}$) $a = 1$ and $b = 0.3338$. The system (20) has only one asymptotically stable limit closed orbits. Its rotation number is 13 (see Fig. 11). It is symmetric to itself about the origin. The period of the limit closed orbit is $T \cong 83.6359$.

$LE_J$:  0.1129, −0.2234, −0.2234  (+,−,−)

$LE_O$: 0.5021, −0.3258, −0.5108  (+,−,−)

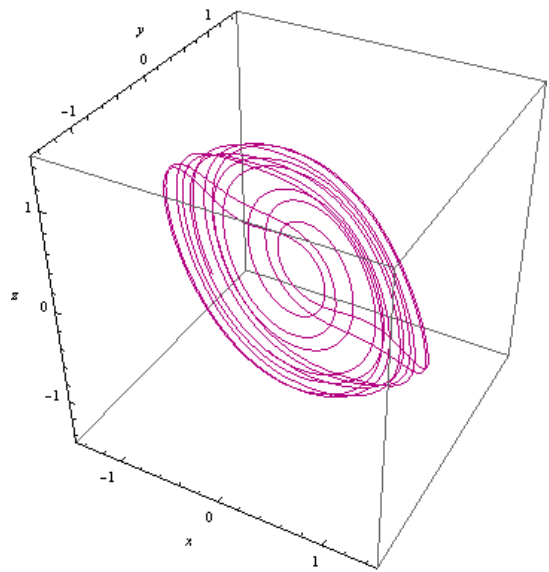

Fig. 11  $a = 1$ and $b = 0.3338$. One limit closed orbit with rotation number 13

From the above numerical results, we see that, for the system (20), $LE_J$ and $LE_O$ are quite different, in the cases ($n_1$), ($n_2$) and ($n_3$), the three exponents of $LE_J$ are all negative, while for $LE_O$, one positive and two negative. In other 8 cases, all of the exponents for the asymptotically limit closed orbits, both $LE_J$ and $LE_O$ have one positive and two negative values, i.e., (+, −, −). So, we have seen that the sign distribution of Lyaponov exponents of a three-dimensional limit closed orbits has the following five possibilities:

(a') (−, −, −)

(b') (0, −, −)

(c') (0, 0, −)

(d') (0, 0, 0)

(e') (+, −, −)

This results is quite different to the proposed criterion for distinguishing attractor given by ref. [1]:

(a) stable fixed point, all the exponents are negative;
(b) limit cycle, an exponent is zero and the remaining ones are all negative;
(c) $k$-dimensional stable torus, the first $k$ LEs vanish and the remaining ones are negative;
(d) for strange attractor generated by a chaotic dynamics at least one exponent is positive.

In addition, for the system (20), the fact that one of $LE_J$ is zero and the others are negative for the asymptotically stable limit cycle happens around the parameter $a=1$ and $b \approx 0.5024$ (see ($n_3$)), and the fact that the limit cycle separated into two happens around the parameter $a=1$ and $b \approx 0.4893$ (see ($n_6$)、($n_7$) and ($n_8$)). These facts show that some exponents equals zero may not be the critical situation that the stability of the limit cycle, or that the number of asymptotically limit cycles changes.

## 4. Conclusions

This paper has shown that there are two kinds of the Lyapunov exponents $LE_J$ and $LE_O$, When they have different values, the second one may lose the basic meaning in dynamical system theory.

The concrete examples have shown that neither of these two exponents could be applied as the numerical characteristic for distinguishing different attractors.

In another paper, the author will give an explanation for why limit closed orbits can still be asymptotically stable when one of the three exponents is positive.

### References

[1] Cencini M. et al., M. *Chaos From Simple models to complex systems*. World Scientific, (2010). ISBN 981-4277-65-7.
[2] Oseledec, V. I. (1968). *A multiplicative ergodic theorem. Lyapunov characteristic numbers for dynamical systems,* Trans. Mosc. Math. Soc. **19**, p. 197.
[3] Keying Guan, *Non-trivial Local Attractors of a Three-dimensional Dynamical System.* arXiv.org, http://arxiv.org/abs/1311.6202(2013).
[4] Keying Guan, Beiye Feng, *Period-doubling Cascades of a Silnikov Equation.* arXiv.org, http://arxiv.org/abs/1312.2043.